\documentclass[11pt]{amsart}

\newtheorem{theorem}{Theorem}

\title{Convexity  and concavity of the ground state energy }
\author{Herbert Koch}
\address{Mathematisches Institut der Universit\"at Bonn \\ Endenicher Allee 60 \\ 53115 Bonn} 
\email{koch@math.uni-bonn.de}  

\begin{document} 

\begin{abstract} This note proves convexity resp. concavity of the ground state 
energy of one dimensional Schr\"odinger operators as a function 
of an endpoint of the interval for convex resp. concave potentials. 
\end{abstract} 

\maketitle

\section{Main result and context}

Let $I=(a,b)\subset \mathbb{R}$ be an open interval, $V \in C(a,b)$ be a convex or concave
potential with $  \liminf_{t \to - \infty} V = \infty$ 
if $a=-\infty$.
   Consider for $t \in (a,b]$ 
the energy 
\[  E_t(u) = \int_a^t u_x^2 + V u^2 dx. \] 
There is a unique positive minimizer $u\in H^1_0(a,t)$ under  
the constraint $\Vert u \Vert_{L^2(a,t)} = 1 $. It satisfies  the Euler-Lagrange  equation
\begin{equation}\label{int}    -u_{xx} + V  u =    \lambda(t)  u \end{equation}   on $(a, t)$ with boundary conditions $u(a)=u(t)=0$ (and obvious modifications if $a=-\infty$). Here $\lambda(t)$ is 
the Lagrangian multiplier, and $\lambda(t)= E_t(u)$. The map $t \to \lambda(t) $ is the main object of interest. 

\begin{theorem} 
The map $(a,b]\ni t \to \lambda(t)$ is twice differentiable, 
strictly decreasing and   $\lim_{t\to a } \lambda(t)=\infty$. 
The map $t\to \lambda(t)$ 
is convex if $V$ is convex, strictly convex if $V$ is convex  and not affine.
If $a=-\infty$ it is  concave if $V$  is concave  and strictly concave if $V$ is concave and not  affine. 
\end{theorem}

 The convexity part follows from a much stronger celebrated result by
 Brascamp and Lieb \cite{MR0450480,BL}.  It is related to a weaker
 statement in Friedland and Hayman \cite{MR0412442} with a computer
 based proof there. These statements found considerable interest and
 use in the context of monotonicity formulas beginning with the
 seminal work of Alt, Caffarelli and Friedman \cite{MR732100}.
 Caffarelli and Kenig \cite{MR1613650} prove a related monotonicity
 formula using the results by Brascamp-Lieb \cite{MR0450480}.  They
 attribute an analytic proof to Beckner, Kenig and Pipher \cite{BKP}
 which the author has never seen.  To the best knowledge of the author
   the concavity statements are new.

 This note has its origin in a seminar of free boundary problems at
 Bonn. It is a pleasure to acknowledge that it would not exist without
 my coorganizer Wenhui Shi. I am grateful to Elliott Lieb for spotting 
an error in the formulation of the main theorem in a previous version.

\section{A short elementary proof} 

\begin{proof} 
Monotonicity and $\lim_{t\to a} \lambda(t) = \infty$ are an immediate 
 consequence of the definition.   We consider the equation \eqref{int} on the interval $(a,t)$ and denote 
  by  $u(x)=u(x,t)$  the unique $L^2$ normalized non negative ground state with ground   state energy $\lambda=\lambda(t) $. Differentiability with respect to 
$x$ and $t$ is an elementary property of ordinary differential equations. 
We argue at a formal level and do not check existence of integrals resp. derivatives below, which follows from  standard arguments. 
We differentiate the equation with respect
  to $t$, denote the derivative of with respect to $t$  by $\dot u$  and obtain
 \begin{equation}    
- \dot u_{xx}  +  V\dot u  -\lambda  \dot u  =  \dot  \lambda    u  
\label{int2} \end{equation}  
with boundary conditions $    \dot u (a)=0$ and $\dot u(t) = -u_x(t)$.  
We multiply \eqref{int2} by $u$, integrate and integrate by parts. Then most terms drop out by \eqref{int}. 
Since $\Vert u \Vert_{L^2}=1$ 
we obtain 
\begin{equation} 
\dot \lambda = \dot u (t) u_x(t) = -   u^2_x(t). 
\label{dynamic} \end{equation} 
 Due to the normalization $\dot u $ is orthogonal to $u$, i.e. $ \int_{a}^t u \dot u  dx = 0$.   
The quotient  $ w = \frac{\dot u }{u}$  satisfies 
\[  w_{xx}   + \frac{u_x}{u} w_x - \frac{u_x^2}{u^2} w  = \dot \lambda  < 0. \] 
In particular $w$ has no non positive local   minimum. 
Since $w \to \infty$ as $x\to t$ 
there can be at most one sign change. Since $\dot u$ is orthogonal to $u$ there 
is exactly one sign change of $\dot u $, lets say at $a<t_0 < t$. Since also 
$\dot u(a)=0$ if $a>-\infty$ we have  $\dot u_x(a)\le 0$ if $a>-\infty$.   
We multiplying \eqref{int}  by $u_x$ and integrate to get 
\begin{equation} \label{derivative}   \dot\lambda = -  u_x(t)^2 =   
   \int_a^t  V'   u^2 dx - u_x^2(a)   \end{equation} 
where we omit the last term here and below if $a=-\infty$.

We differentiate \eqref{derivative} with respect to $t$ and 
 use the orthogonality  $\int_a^t u \dot udx = 0 $ to obtain a partly implicit  formula for the second derivative of $\lambda$ with respect to $t$,   
\[
\begin{split} 
 \ddot\lambda = &   2 \int_a^t  (V'(x)-V'(t_0)) u  \dot u  dx - 2 u_x(a) \dot u_x(a) 
\\ = &   2\int_a^t   (V'(x)-V'(t_0))  w u^2 dx  -2 u_x(a)\dot u_x(a). 
\end{split} 
\] 
Recall that $u_x(a)>0$ and $\dot u_x(a)\le 0 $ and hence the second term on the right hand side is nonnegative. By the choice of $t_0$  the first term is nonnegative if $V$ is convex, nonpositive if it is concave, 
positive if $V$ is convex and non affine, and negative if $V$ is concave and not affine.  
Thus $t\to \lambda$ is convex if $V$ is convex, it satisfies $\ddot \lambda >0$ 
if $V$ is convex and not affine (i.e. $V'$ is not constant), if $a=-\infty$ 
it is concave if $V$ is concave and $\ddot \lambda < 0$ if $V$ is concave 
and  not affine. 
\end{proof}

\bibliography{seminar}{} 
\bibliographystyle{plain} 

\end{document}